\documentclass[11pt]{article}
\usepackage{graphicx}
\usepackage[T1]{fontenc}
\usepackage[utf8]{inputenc}
\usepackage{mathtools}
\usepackage{amsmath}
\usepackage{amssymb}
\usepackage[margin=2.5cm]{geometry}
\usepackage{amsfonts}
\usepackage{amsmath,amssymb,amsthm}
\RequirePackage{luatex85}
\usepackage{enumitem}
\title{Solutions of the Exponential Equation $7x^2 + 59y^2 = 3^m$ \\
A simple Algorithm producing all the primitive solutions}
\author{Roy Barbara}
\date{}
\begin{document}
\maketitle
\begin{center}
\footnotesize{A}\tiny{BSTRACT.}\hspace{1.3mm}\footnotesize{We provide an \textit{elementary} algorithm, with no use of radicals or complex numbers, \\ and with \textit{elementary} proof, that generates \ \textit{all} \ the \ (infinitely many) \ primitive and\hspace{1.5mm}positive\\  solutions of the exponential\hspace{0.6mm} \hspace{0.6mm}diophantine equation \ $7x^2 + 59y^2 = 3^m.$ \ \ \ \ \  \ \  \ \ \  \ \ \ \ \ \  \ \  \ \ \  \ \ \ \ \ \  \ \  \ \ \hspace{0.03mm}}
\end{center}
${}$ \\ \\
{\Large \bf{1. Introduction}} \\ \\
A diophantine equation is a polynomial equation in two or more variables (unknowns), representing integers (or sometimes rational numbers). Examples are: $x^2 + y^2 = 2z^6$, $x^3 + y^4 = z^5$, the Pell's equation $x^2 - dy^2 = \pm 1$, \textit{etc.} If one (or more) of the variables occurs as an exponent, the equation is called an exponential diophantine equation. Examples are: $2x^2 + 3y^2 = 5^z$, $3^x + 4^y = z^2$, the Ramanujan-Nagell equation $2^n - 7 = x^2$, \textit{etc.} A purely exponential (diophantine) equation is one in which \textit{all} the variables occur as exponents, as for example $2^x + 3^y = 5^z$. There is no general method to solve an exponential diophantine equation. Starting with purely exponential (ternary) equations as $a^x + b^y = c^z$, we observe that such equations have in general a (small) \textit{finite} number of solutions: In $[1]$, it is proven that the equation $3^x + 4^y = 5^z$ has the single positive integral solution $(x,y,z) = (2,2,2)$. (this result has been generalized to other pythagorian triplets). Scott proved that the equation $3^x + 13^y = 2^z$ has exactly \textit{two} positive solutions, and that the equation $3^x +5^y = 2^z$ has exactly \textit{three} positive solutions. (see $[2]$). Further, it is conjectured in $[2]$ that the equation $a^x + b^y = c^z$ has at most \textit{one} solution with $x,y,z \ge 2$. (this conjecture has been revisited). We move next to the equations in \textit{two} exponent variables, \ a typical form being \ \ $a^x + b^y = z^2$. \textit{Numerous} results can be found: For example, in $[3]$, it is proven that the equation $8^x + 19^y = z^2$ has a single nonnegative integer solution, and in $[4]$, it is proven that the equation $2^x + 17^y = z^2$ has exactly \textit{five} positive integer solutions. As far as one can see, such equations again have in general a (small) \textit{finite} number of solutions. Finally, we focus on exponential diophantine equations in just \textit{one} exponent variable. A typical form is $ax^2 + by^2 = \lambda k^z$ ($\lambda = 1,2,4$). Such equations \textit{would} be easy to solve \textit{if} the integers represented by the quadratic form $ax^2 + by^2$ \textit{formed} a multiplicative semi-group (we disregard the trivial case where $a=1$ or $b=1$, with $\lambda = 1$). Several results can be found (for ex. see $[5]$). Elaborated recent articles are also available. \ However, \textit{totally absent}${}$ is a simple and practical algorithm, with elementary proof, that generates all the (infinitely many) primitive solutions of a \textit{given} exponential equation as for example $5x^2 + 7y^2 = 3^z$, or $3x^2 + 5y^2 = 17^z$, etc. In $2018$, an open problem appeared in different Fb groups, including Terence Tao's fan club, asking for the solutions in integers $x,y,m$ $(m>0)$ of the exponential equation $7x^2 + 59y^2 = 3^m$, with $gcd(x,y) = 1$. We have several motivations towards this problem: First, this uncommon problem seems to be attractive. The choice of the primes $7,59$ and $3$ is not irrelevant. We note for example that the equation $7x^2 + 59y^2 = 3^m$ has \textit{no} integral solution with $z$ \textit{even}, and then that the Gauss' composition does not apply to the set of integers represented by the quadratic form $7x^2 + 59y^2$.
\pagebreak

Since then, one wonders \textit{how} and \textit{where} to start in order to find the first primitive solutions (that anyway grow exponentially). Secondly, in contrast with other kinds of exponential equations, it turns out that the equation $7x^2 + 59y^2 = 3^m$ possesses infinitely many primitive solutions; that encourages us to search for a simple algorithm, using elementary tools, which will generate all these primitive solutions. Thirdly, our approach, using new ideas and an appropriate "\textit{finite descent}", can be generalized to provide a complete solution, in an elegant and \textit{effective} way, to a large (infinite) class of exponential diophantine equations of the form $ax^2 + by^2 = \lambda k^z$, $a,b,k > 1$.  \\ \\ \\
{\Large \bf{2. The Problem }}\\ \\
In this article, we consider the exponential diophantine equation
\begin{equation} 7x^2 + 59y^2 = 3^m \tag{1}
\end{equation}
The variables $x,y$ and $m$ denote non-zero integers, with $m>0$. \\ A solution $(x,y,m)$ of $(1)$ is said to be "positive" if $x$ and $y$ are positive.  To avoid trivialities, \\ we only need to focus on primitive solutions (i.e. solutions where $x$ and $y$ are coprime).\\ \\
\textit{Notation}:\\
A "pp-solution" of $(1)$ stands for: a \textit{primitive} and \textit{positive} solution of $(1)$.\\ \\
Note that if we know all the pp-solutions of $(1)$, hence all the primitive solutions of $(1)$, we know all the solutions of $(1)$ (just observe that, in any solution of $(1)$, the gcd of $x$ and $y$ must be a power of $3$). \\ \\
\textit{Definition}:\\
We say that a positive integer $m_0$ is "suitable" if there is a primitive solution of $(1)$ of the form $(x,y,m_0)$.\\ \\
The least suitable integer is $5$, corresponding to the single pp-solution of $(1)$, $S_1 = (1,2,5)$. One may check that the least suitable integer greater than $5$, is $15$, corresponding to the single pp-solution \\ of $(1)$, $S_2 = (701, 430, 15)$.\\\\
Natural questions arise: \begin{itemize}
\item[-] Are there finitely or infinitely many primitive solutions of $(1)$?
\item[-] Which positive integers are suitable?
\item[-] If $m_0$ is suitable, is there a \textit{unique} pp-solution of $(1)$, of the form $(x,y,m_0)$?
\item[-] Is it possible to produce \textit{all} the primitive solutions of $(1)$ by a simple algorithm?
\item[-] Is there any solution of $(1)$ with $x$ even?
\end{itemize}
We will answer all these questions by proving the following: \\ \\
There are \textit{infinitely many} primitive solutions of $(1)$, and we are able to \textit{determine all} these in an elementary \textit{effective way}. It turns out that the suitable positive integers $m$, \textit{all} odd, form precisely the arithmetic progression:
$$m = 10k + 5, \ \ \ k = 0,1,2,3, \ldots$$
and that, for each given $m_0 = 10k_0 + 5$, there is \, \textit{exactly one} \, pp-solution of the form $(x,y,m_0)$. \\ The first pp-solutions of $(1)$ are:
$$(1,\ 2,\ 5), \ (701, \ 430,\ 15), \ (262009, \ 78842, \ 25), \ (78606773, \ 10718566, \ 35)\ .$$
At the end, we deduce that there is no solution of $(1)$ with $x$ even. \\ \\
\textit{Notation}:\\
Let $S = (x,y,m)$ and $S' = (x',y',m')$ be two pp-solutions of $(1)$. We say that $S$ is "smaller" than $S'$,\, and we write $S < S'$, \, to mean that $m<m'$. \\ \\ \\
{\Large \bf{3. The Results}} \\ \\
\textit{Notation}:\\
If $p$ and $q$ are non-zero integers, we define the integers:
\begin{align*}
A(p,q) & = 59p^2 - 236pq - 7q^2  \\
B(p,q) & = -118p^2 -14pq + 14q^2 \\
C(p,q) & =  9(59p^2 + 7q^2)
\end{align*}
{\large \textbf{Theorem 1:}}\\
Let $(q,p,\omega)$ be a pp-solution of $(1)$, so that ${}$ $7q^2 + 59p^2 = 3^{\omega}$. ${}$ Clearly, $3$ doesn't divide $pq$, so that \\ $pq \equiv \pm 1 \pmod3$. \textit{In the formulas below} (replacing $p$ by $-p$ \textit{if necessary}), we may always assume that $pq \equiv -1 \pmod3$. Set $$x = 3^{-5} |A(p,q)|,\ \  y = 3^{-5} |B(p,q)|, \ \  x' = |A(-p,q)|, \ \ y' = |B(-p,q)|.$$
Then, \\$(x,y,2\omega -5)$ is a pp-solution of $(1)$, that we call the first successor of $(q,p,\omega)$, \\and\\
$(x',y',2\omega+5)$ is a pp-solution of $(1)$, that we call the second successor of $(q,p,\omega).$ \\ \\
{\large \textbf{Remark 1:}}\\
The first successor of $S_1 = (1,2,5)$ is $S_1$ itself and the second successor of $S_1$ is $S_2 = (701, 430, 15)$. Now, let $S = (q,p,\omega)$ be a pp-solution of $(1)$, with $S\neq S_1$ (so $\omega \ge 15$). As quickly seen, if $S'$ and $S''$ denote respectively the first and second successor of $S$, we then have \, $S < S' < S''$. \\ \\
{\large \textbf{Corollary 1:}}\\
By starting with $S_1$ and $S_2$, then by taking the successors of $S_2$, say $S_3 < S_4$, and then by taking the successors of $S_3$, say $S_5 < S_6$, and the successors of $S_4$, say $S_7 < S_8$, and so on, we obtain an infinite binary tree of pp-solutions of $(1)$. The reader can easily check that the suitable $m$'s obtained this way, form precisely the arithmetic progression $$5,\ 15,\ 25,\ 35,\ 45,\ 55,\ 65,\ 75, \ \ldots$$
In particular, equation $(1)$ possesses infinitely many primitive solutions. \\ \\ \\
{\large \textbf{Theorem 2:}}\\
Every pp-solution of $(1)$ belongs to the binary tree described above. In other words, the previous algorithm generates \textit{all} the primitive and positive solutions of $(1)$. \\ \\
${}$ \; \; \; \; $\bullet$ Finally, we observe that in every solution $(x,y,m)$ of $(1)$, \, $x$ is odd (and $y$ is even). \\ \\ \\
{\Large \bf{4. Proof of theorem 1}}\\ \\
\textit{Notation}:\\ $\mathbb Z^{\ast}$ will denote the set of non-zero integers. \\ \\
{\large \textbf{Lemma 1:}}\\
Let $p,q \in \mathbb Z^{\ast}$. Then,
\begin{enumerate}
\item[(i)] We have the identity:
$$7(59p^2 -236pq -7q^2)^2 + 59(-118p^2 - 14pq + 14q^2)^2 = 3[9(59p^2 +7q^2)]^2$$
Stated in a compact way: if $A = A(p,q), B= B(p,q)$ and $C = C(p,q)$, \, then,
$$7A^2 + 59B^2 = 3C^2$$
\item[(ii)] If we suppose further that $7q^2 + 59p^2 = 3^{\omega}$, then
$$7A^2 + 59B^2 = 3^{2\omega + 5}$$
\end{enumerate}
\begin{proof}
\hfill
\begin{enumerate}
\item[(i)] The proof of the identity is straightforward.
\item[(ii)] $7q^2 + 59p^2 = 3^{\omega}$ \, yields \, $C = C(p,q) = 9(7q^2 + 59p^2) = 3^{\omega+2}$. \, Hence, $(i)$ yields \\$7A^2 + 59B^2 = 3C^2 = 3(3^{\omega+2})^2 = 3^{2\omega + 5}$.
\end{enumerate}
\end{proof}
{\large \textbf{Remark 2:}}\\
Let $p,q \in \mathbb Z^{\ast}$. Then, $A(p,q) \neq 0$ \, \textit{and} \, $B(p,q) \neq 0$.\\
Indeed: set $u = A(p,q)$, \,$v = B(p,q)$ \, and \, $w = C(p,q)$. \, Then, $w = 9(59p^2 + 7q^2) > 0$,\, and, by \\part $(i)$ of lemma $1$, we have $7u^2 + 59v^2 = 3w^2$. Assuming $u=0$ would lead to the contradiction that $\sqrt{\frac{3}{59}}$ is rational, and assuming $v = 0$ would lead to the contradiction that $\sqrt{\frac{3}{7}}$ is rational. \\ \\
{\large \textbf{Lemma 2:}}\\
Let $u,v \in \mathbb Z^{\ast}$ and let $m$ be an integer, $m \ge 15$. Suppose that
\begin{enumerate}
\item[(i)] $7u^2 + 59v^2 = 3^m$, where $3 \mid u$.
\item[(ii)] $2u + v = 3^5 \lambda$ \,\, for some $\lambda \in \mathbb Z^{\ast}$, \, where $3 \nmid \lambda$.
\end{enumerate}
Then, \, $\gcd(u,v) = 3^5$.
\begin{proof}\hfill \\
Since $3\mid u$, then by $(i)$, we see that $3 \mid v$. Relation $(i)$ also shows that any common prime factor of $u$ and $v$ must be $3$. Therefore,
$$\gcd(u,v) = 3^r, \; \; r \ge 1$$
If we assume that $r \ge 6$, then $3^6$ would divide $2u+v$, and hence by $(ii)$, $3^6$ would divide $3^5 \lambda$, so $3$ would divide $\lambda$, a contradiction. We conclude that $$\gcd(u,v) = 3^r, \; \; 1 \le r \le 5$$
Now, we claim that $r = 5$: for the purpose of contradiction, suppose that $1 \le r \le 4$. Under this assumption, $3^5$ cannot divide $u$: otherwise by $(i)$ and as $m\ge 15$, ${}$ $3^5$ would divide $v$, hence $3^5$ would divide $\gcd(u,v)$. \, We then set
\begin{equation}
u = 3^r \theta, \; \; \; \; \text{where}\; 1 \le r \le 4 \; \; \text{and} \; \; 3 \nmid \theta \tag{j}
\end{equation}
Using $(i)$ and $(ii)$ , we may write
$$7u^2 + 59(3^5 \lambda - 2u)^2 = 3^m$$
Expanding yields
$$7u^2 + 3^{10} \cdot 59 \lambda^2 + 236u^2 - 3^5\cdot 236 \lambda u = 3^m$$
That is,
$$3^5 u^2 + 3^{10} \cdot 59 \lambda^2 - 3^5 \cdot 236 \lambda u = 3^m$$
Replacing $u$ by $3^r \theta$ (given by $(j)$) gives
$$3^{5 + 2r} \cdot \theta^2 + 3^{10} \cdot 59  \lambda^2 - 3^{5+r} \cdot 236 \lambda \theta = 3^m$$
Dividing by $3^{5+r}$ yields
\begin{equation} 3^r \theta^2 + 3^{5-r} \cdot 59 \lambda^2 - 236 \lambda \theta = 3^{m-5-r} \tag{$\ell$}
\end{equation}
Since $1\le r\le 4$ \, and \, $m\ge 15$, \, all the terms in $(\ell)$ are divisible by $3$, \textit{except} the term $-236 \lambda \theta$. \\ We obtain a contradiction.
\end{proof}
\hfill \\
{\large \textbf{Corollary 2:}}\\
Let $(q,p,\omega)$ be a primitive solution of $(1)$, with $pq \equiv -1 \pmod3$. Set $u = A(p,q)$ and $v = B(p,q)$. Then, $\gcd(u,v) = 3^5$.
\begin{proof} \hfill \\
We have $7q^2 + 59p^2 = 3^{\omega}$ ($\omega \ge 5$). Set $m = 2\omega + 5$, so $m\ge 15$. By lemma $1$ $(ii)$, we have $7u^2 + 59v^2 = 3^m$. Using $p^2 \equiv q^2 \equiv 1 \pmod3$ and $pq \equiv -1 \pmod3$, we may write $u = 59p^2 - 236pq -7q^2 \equiv 59 + 236-7 = 288  \equiv 0 \pmod3$. ${}$ Hence, $3$ divides $u$. ${}$ On the other hand, we have
\begin{align*}
2u + v &= 2(59p^2 - 236pq - 7q^2) + (-118p^2 -14pq + 14q^2) \\
&= -486 pq = 3^5(-2pq), \; \; \text{where} \; 3 \nmid -2pq
\end{align*}
By lemma $2$, we obtain that $\gcd(u,v) = 3^5$.
\end{proof} \hfill \\
\textit{Proof of theorem 1}:\\
Let $(q,p,\omega)$ be a pp-solution of $(1)$. As pointed out in theorem $1$, we may assume that $pq \equiv -1 \pmod 3$. Set $u = A(p,q)$ and $v = B(p,q)$. By corollary $2$, we have $\gcd(u,v) = 3^5$. Set $x = 3^{-5} |u|$ and $y = 3^{-5} |v|$. Thus, $x$ and $y$ are \textit{coprime positive} integers.\\
Since by hypothesis $7q^2 + 59p^2 = 3^{\omega}$, lemma $1$ $(ii)$ provides:
$$7u^2 + 59v^2 = 3^{2\omega + 5}$$
Dividing by $3^{10}$ yields
$$7(3^{-5} u)^2 + 59(3^{-5} v)^2 = 3^{2\omega - 5}$$
Hence,
$$7x^2 + 59y^2 = 3^{2\omega - 5}$$
As $x > 0$, \, $y>0$ \, and \, $\gcd(x,y) = 1$, \, we see that $(x,y,2w - 5)$ is a pp-solution of $(1)$.\\\\
$\bullet$ Next, set $x' = |A(-p,q)|$ and $y' = |B(-p,q)|$.\\
With $pq \equiv -1 \pmod3$, we may write \\
$A(-p,q) = 59p^2 + 236 pq - 7q^2 \equiv 59 -236 - 7 = -181 \equiv -1 \pmod3$.\\
Hence, $3 \nmid x'$. Now, lemma $1$ $(ii)$ applied to $-p$ and $q$ yields:
$$7[A(-p,q)]^2 + 59[B(-p,q)]^2 = 3^{2\omega + 5}$$
Finally, we have $7x'^2 + 59y'^2 = 3^{2\omega + 5}$, where $3 \nmid x'$.\\
Hence clearly, $(x',y', 2\omega +5)$ is a pp-solution of $(1)$. \hfill $\square$ \\ \\ \\
{\large \textbf{Remark 3:}}\\
In the light of what precedes and particularly of lemma $1$ $(ii)$, we end this section by providing a simple and useful criterion to "\textit{recognize}" the successor, that we will use in section $6$. \\ \\
$\bullet$ Given that $(x,y,m)$ and $(q,p,\omega)$ are known to be two \textit{primitive} and \textit{positive} solutions of $(1)$  (here, we need neither to suppose that $pq \equiv -1 \pmod3$, nor to assume any relation between the exponents $m$ and $\omega$), then, the reader can check the following:
\begin{enumerate}
\item[(i)] If we can just show that $x = |A(\pm p, q)|$ and $y = |B(\pm p, q)|$ (the "plus" is to be taken with the "plus", and the "minus" with the "minus"), then, $(x,y,m)$ is the second successor of $(q,p,\omega)$.
\item[(ii)] If we can just show that $x = 3^{-5} |A(\pm p,q)|$ and $y = 3^{-5} |B(\pm p,q)|$ (the "plus" is to be taken with the "plus", and the "minus" with the "minus"), then, $(x,y,m)$ is the first successor of $(q,p,\omega)$. \\
\end{enumerate}
{\Large \bf{5. Preliminaries for theorem 2:}} \\ \\
\textit{Notation}: \\ As usual, $\left( \frac{x}{y} \right)$ denotes Legendre's symbol. \\ \\
Recall Legendre's theorem for ternary quadratic forms (see for ex. $[6]$):\\ \\
\textbf{Legendre's Theorem:}\\
Let $a,b,c$ be positive square-free integers, and pairwise coprime. Then, the equation $ax^2 + by^2 = cz^2$ has a non-trivial integer solution \; \textit{if and only if} \; $\left( \frac{-ab}{c} \right) = \left( \frac{bc}{a} \right) = \left( \frac{ca}{b} \right) = +1$. \\ \\
{\large \textbf{Lemma 3:}}\\
The equation $7x^2 + 59y^2 = z^2$ has \textit{no} non-trivial integer solution.
\begin{proof} \hfill \\
Since $\left(\frac{59 \times 1}{7} \right) = \left( \frac{3}{7} \right) = -1$, the result follows from Legendre's theorem.
\end{proof}
${}$ \, \, \, \, $\bullet$ Since $3^m$ has the form $z^2$ when $m$ is even, we obtain: \\ \\
{\large \textbf{Corollary 3:}} \\
The equation $7x^2 + 59y^2 = 3^m$ has \textit{no} integer solution with $m$ even. \\ \\
{\large \textbf{Lemma 4:}} \\
The equation $7x^2 + 59y^2 = 2z^2$ has \textit{no} non-trivial integer solution.
\begin{proof} \hfill \\
Since $\left( \frac{59 \times 2}{7} \right) = \left( \frac{-1}{7} \right) = -1$, the result follows from Legendre's theorem.
\end{proof} \hfill \\
{\large \textbf{Lemma 5:}} \\
The equation $7x^2 + 59y^2 = 6z^2$ has \textit{no} non-trivial integer solution.
\begin{proof} \hfill \\
We show that $\left( \frac{6\times 7}{59} \right) = -1$, and the result follows from Legendre's theorem.\\We have $\left( \frac{6\times 7}{59} \right) = \left( \frac{2}{59} \right) \left( \frac{3}{59} \right) \left(\frac{7}{59} \right)$. \\
Now, using (Gauss') law of quadratic reciprocity, etc. (see for ex. $[7]$ or $[8]$) we may write:
\begin{align*}
\left( \frac{2}{59} \right) &= (-1)^{\frac{59^2 -1}{8}} = (-1)^{435} = -1. \\
\left(\frac{3}{59}\right) &= \left(\frac{59}{3} \right) (-1)^{\frac{3-1}{2} \cdot \frac{59 -1}{2}} = \left( \frac{2}{3} \right)(-1)^{29} = -\left(\frac{2}{3}\right) = -(-1) = +1. \\
\left(\frac{7}{59} \right) &= \left(\frac{59}{7}\right) (-1)^{\frac{7-1}{2} \cdot \frac{59-1}{2}} = \left(\frac{3}{7} \right) (-1)^{87} = -\left( \frac{3}{7} \right) = -(-1) = +1.
\end{align*}
Hence, $\left( \frac{6\times 7}{59} \right) = (-1)(+1)(+1) = -1$.
\end{proof}
\hfill \\
${}$ \; \; \; \; $\bullet$ Since $2 \cdot 3^m$ has the form ${}$ $2z^2$ ${}$ or ${}$ $6z^2$ ${}$ according to whether $m$ is even or odd, we obtain \\ (by lemmas $4$ and $5$) the following: \\ \\
{\large \textbf{Corollary 4:}} \\
The equation $7x^2 + 59y^2 = 2\cdot 3^m$ ($m > 0$) has \textit{no} integer solution. \\ \\
{\large \textbf{Remark 4:}}\\
Let $p,q \in \mathbb Z^{\ast}$. Set $A = A(p,q)$, $B = B(p,q)$ and $C = C(p,q)$.\\
Let $\delta$ denote the positive gcd of $A$ and $B$. According to lemma $1$ $(i)$, we have: $7A^2 + 59B^2 = 3C^2$.
In virtue of this identity, it should be clear that:
$$\delta = \gcd(A,B) = \gcd(A,C) = \gcd(B,C). \\\\\\$$\\\\
{\large \textbf{Proposition 1:}}\\
Let $p,q \in \mathbb Z^{\ast}$, $p$ and $q$ coprime. Set $A = A(p,q)$ and $B = B(p,q)$. Let $\delta$ denote the positive gcd of $A$ and $B$. Then,
\begin{enumerate}
\item[(i)] $\delta$ divides $2\cdot 3^5 \cdot 7 \cdot 59$.
\item[(ii)] We have: \; $[ \ 7 \mid \delta \iff 7 \mid p \ ]$\;\; and\;\; $[ \ 59 \mid \delta \iff 59 \mid q \ ]$. \\
\end{enumerate}
\textit{Proof}:
\begin{enumerate}
\item[(i)] If $\delta = 1$, there is nothing to prove. We then assume that $\delta > 1$. We show that every \textit{primary} factor of $\delta$ must have one of the following forms:
$$2^1, \; \; \; \; 3^s \; (1 \le s\le 5), \; \; \; \; 7^1 \; \;\;\; \text{or} \;\; \; \; 59^1$$
Indeed, let $\pi^s$, $s \ge 1$, be a \textit{primary} factor of $\delta$. Then, $\pi^s$ divides $A$ and $B$. Hence, $\pi^s$ divides $2A + B = -2 \cdot 3^5 \cdot pq$. We consider two cases:\\ \\
\textit{case 1}: $\pi \nmid pq$:\\
Then, $\pi^s$ is coprime to $pq$. By Gauss' theorem, $\pi^s$ divides $2 \cdot 3^5$. Hence, either $\pi = 2$ and $s = 1$,\;  or, \; $\pi = 3$ \;  and \; $1 \le s \le 5$.\\ \\
\textit{case 2}: $\pi \mid pq$:\\
Since further $p$ and $q$ are coprime, there are exactly two possibilities:\\
sub-case 1: $\pi \mid p$ and $\pi \nmid q$:\\
From $\pi \mid p$ and $\pi \mid A = p(59p^2 - 236q) - 7q^2$, we deduce that $\pi \mid 7q^2$. As further $\pi \nmid q$, then, $\pi = 7$. In particular, $7 \nmid q$. Now, $\pi^s = 7^s$ divides $2A + B = -2\cdot 3^5 \cdot pq$, where $7^s$ is coprime to $2 \cdot 3^5 \cdot q$. By Gauss' theorem, $7^s$ divides $p$. Finally, from $7^s \mid p$ and $7^s \mid A = p(59 p-236q)-7q^2$, we deduce that $7^s \mid 7q^2$. Since further $7^s$ is coprime to $q^2$, we get, $7^s \mid 7$. ${}$ Hence $s = 1$ (so $\pi^s = 7^1$).\\
sub-case 2: $\pi \mid q$ and $\pi \nmid p$:\\
From $\pi \mid q$ \; and \; $\pi \mid A = q(-236p - 7q) + 59p^2$, we deduce that $\pi \mid 59p^2$. As further $\pi \nmid p$, then, $\pi = 59$. In particular, $59\nmid p$. Now, $\pi^s = 59^s$ divides $2A +B = -2\cdot 3 \cdot pq$, where $59^s$ is coprime to $2\cdot 3^5 \cdot p$. \; By Gauss' theorem, $59^s$ divides $q$. \; Finally,\hspace{2mm}from \; $59^s \mid q$ \; and \\ $59^s \mid A = q(-236p -7q^2) + 59p^2$, \; we deduce that $59^s \mid 59p^2$.\\
Since further $59^s$ is coprime to $p^2$, we get ${}$ $59^s \mid 59$. \; Hence $s = 1$ (so $\pi^s = 59^1$).
\item[(ii)] \begin{enumerate}
\item[$\bullet$] If $7 \mid \delta$, then $7\mid B = 7(-2pq + 2q^2) - 118p^2$. Hence, $7\mid 118p^2$. As further $7\nmid 118$, then, \\ $7$ must divide $p$.\\
Conversely, if $7 \mid p$, then, $7\mid A = p(59p - 236q) -7q^2$\; and \; $7\mid B = p(-118p - 14q) + 7(2q^2)$. Hence, $7 \mid \delta$.
\item[$\bullet$] Similarly, if $59 \mid \delta$, then $59 \mid A = 59(p^2 -4pq) - 7q^2$. Hence, $59 \mid 7q^2$. Then, clearly, $59 \mid q$. Conversely, if $59 \mid q$, then $59 \mid A = 59p^2 + q(-236 p - 7q)$ and $59 \mid B = 59(-2p^2) + q(-14p + 14q)$. Hence, $59 \mid \delta$.
\end{enumerate} \hfill $\square$
\end{enumerate}
\hfill \\
Next, we give a parametrization of the rational points on the ellipse $7X^2 + 59Y^2 = 3$: Given $p,q \in \mathbb Z^{\ast}$, set $A  = A(p,q)$, $B = B(p,q)$ and $C = C(p,q)$. In virtue of lemma $1$ $(i)$, we have $7A^2 + 59B^2 = 3C^2$. Since $C > 0$, we obtain $7\left( \frac{A}{C} \right)^2 + 59 \left( \frac{B}{C}\right)^2 = 3$. Hence, $(X,Y) = \left(\frac{A}{C}, \frac{B}{C} \right)$ is a rational point on the ellipse $7X^2 + 59Y^2 = 3$. \\\\
Now, we focus with more attention and precision on the converse: \\

{\large \textbf{Proposition 2:}} \\
Let $(X,Y)$ be a rational point on the ellipse $7X^2 + 59Y^2 = 3$, such that $X \neq \frac19$ and $Y \neq \frac29$. Then, we have
$$X = \frac{A(p,q)}{C(p,q)} \; \; \;\; \; \; \; \; \text{and} \; \; \; \; \;\;\;\; Y = \frac{B(p,q)}{C(p,q)}$$
for some $p,q \in \mathbb Z^{\ast}$, \; $q>0$\;, \; $p$ and $q$ coprime, and where the fraction (in lowest terms) $\frac{p}{q}$ is \textit{precisely} equal to the non-zero rational number $\theta = \frac{9Y-2}{9X-1}$. \\
\begin{proof} \hfill \\
Let $(X,Y)$ be a rational point on the ellipse $7X^2 + 59Y^2 = 3$, with $X \neq \frac19$ and $Y \neq \frac29$. Define the non-zero rational number $\theta = \frac{9Y-2}{9X-1}$. \; Solving for $Y$, we get
\begin{equation}
Y = \; \theta X + \frac{2-\theta}{9} \tag{4.1}
\end{equation}
Using $(4.1)$, we then have, \;\; $7X^2 + 59\left(\theta X +\frac{2-\theta}{9}\right)^2 = 3$.
\; Multiplying by $81$ yields
\begin{equation}
7(9X)^2 + 59(9\theta X + 2 - \theta)^2 = 243 \tag{4.2}
\end{equation}
Expanding $(4.2)$ and after a little algebra, we find
$$81(59\theta^2 + 7)X^2 + 1062 \theta(2 - \theta) X + (59 \theta^2 - 236 \theta -7) = 0$$
The reduced discriminant of this trinomial in $X$ is $\Delta = 81(118 \theta + 7)^2$. Hence, $\sqrt{\Delta} = 9|118 \theta + 7|$, so,
$$\{\pm \Delta \} = \{ \pm 9(118 \theta + 7) \}.$$
According to the hypothesis, the root ${}$ $X = \frac19$ \; has to be rejected.\hspace{3mm}Therefore, we must have
\begin{equation}
X = \frac{59 \theta^2 - 236 \theta - 7}{9(59 \theta^2 + 7)} \tag{4.3}
\end{equation}
Using $(4.1)$ and $(4.3)$, we obtain
\begin{equation}
Y = \frac{-118 \theta^2 - 14 \theta + 14}{9(59\theta^2 + 7)} \tag{4.4}
\end{equation}
Finally, let us write $\theta$ in lowest terms as $\theta = \frac{p}{q}$, where $p,q \in \mathbb Z^{\ast}$, $q>0$ ($p$ and $q$ coprime). \\
Replacing $\theta$ by $\frac{p}{q}$ in $(4.3)$ and $(4.4)$ yields
$$X = \frac{59p^2 - 236 pq - 7q^2}{9(59 p^2 + 7q^2)} = \frac{A(p,q)}{C(p,q)}$$
and
$$Y = \frac{-118p^2 - 14pq + 14q^2}{9(59p^2 + 7q^2)} = \frac{B(p,q)}{C(p,q)}.$$
\end{proof}
$\bullet$ Finally, we study the diophantine equation
\begin{equation}
7x^2 + 59y^2 = 3z^2 \tag{2}
\end{equation}
where the variables $x,y,z$ are non-zero integers, and where $z>0$. A solution of $(2)$ is said to be \textit{positive} when furthermore $x$ and $y$ are positive. A solution $(x,y,z)$ of $(2)$ is called \textit{primitive} when $x$ and $y$ are coprime (equivalently, when $x,y,z$ are pairwise coprime). \\\\
\textit{Notation}:\\
A "pp-solution" of $(2)$ stands for: a \textit{primitive} and \textit{positive} solution of $(2)$.\\ \\
Note first that if $p,q \in \mathbb Z^{\ast}$, and if $x = A(p,q)$, $y = B(p,q)$ and $z = C(p,q)$, then, due to lemma $1(i)$, $(x,y,z)$ is a solution of $(2)$. However, even if $p$ and $q$ are coprime, this solution is not necessarily primitive. \\
We are particularly interested in the converse, that we study with more attention and precision: \\ \\
\textit{Definition}:\\
Let $S = (x,y,z)$ be a primitive solution of $(2)$, with $z>0$, such that,
$$(x,y,z)\;\; \neq \;\; (1,2,9), \;\; (1,-2,9), \;\; (-1,2,9).$$
We define the \textit{incidence} of $S$, that we denote by $I(S)$, as the non-zero rational number
$$I(S) = \frac{9y -2z}{9x -z}$$ \hfill \\
\textit{Clarification}:\\
We clarify why $I(S)$ is well-defined and why $I(S) \neq 0$:\\
If we suppose that $9x-z = 0$, then we get $ \frac{x}{z} = \frac19$. Since the fractions $\frac{x}{z}$ and $\frac19$ are \textit{both} in lowest terms, and since $z>0$, we obtain $z = 9$ and $x = 1$. From this and $7x^2 + 59y^2 = 3z^2$, we get $y = \pm 2$, so $(x,y,z) = (1,\pm 2, 9)$, a contradiction.\\
On the other hand, if we suppose that $9y - 2z = 0$, then we get $\frac{y}{z} = \frac29$. Since \textit{both} fractions $\frac{y}{z}$ and $\frac29$ are in lowest terms, and since $z>0$, we obtain $z = 9$ and $y = 2$. From this and $7x^2 + 59y^2 = 3z^2$, we get $x = \pm 1$, so $(x,y,z) = (\pm 1,2,9)$, a contradiction. \\ \\
{\large \textbf{Proposition 3:}}\\
Let $S=(x,y,z)$ be a primitive solution of $(2)$, with $z>0$, such that
$$(x,y,z) \;\; \neq \;\; (1,2,9), \; \; (1,-2,9) \;\; (-1,2,9).$$
Write the incidence $I(S)$ in lowest terms as $I(S) = \frac{p}{q}$, with $q>0$. Set $A = A(p,q)$, $B = B(p,q)$, $C = C(p,q)$, and let $\delta$ be the positive gcd of $A$ and $B$. Then, we have \textit{precisely}:
$$x = \frac{A}{\delta}, \;\;\;\; y = \frac{B}{\delta} \;\;\;\; \text{and} \;\;\;\; z = \frac{C}{\delta}$$
where $\delta$ \textit{divides} $2\cdot 3^5 \cdot 7 \cdot 59$.\\
Furthermore, we have \; \; [\;$7\mid \delta \iff 7 \mid p$\;] \; \; and \; \; [\;$59 \mid \delta \iff 59 \mid q$\;].
\begin{proof} \hfill \\
Set $X = \frac{x}{z}$ and $Y = \frac{y}{z}$. We have $7x^2 + 59y^2 = 3z^2$. Dividing by $z^2 \neq 0$ yields $7X^2 + 59Y^2 = 3$. Since $I(S)$ is defined and is non-zero, we have $9x - z \neq 0$ and $9y - 2z \neq 0$. Hence, $X \neq \frac19$ and $Y \neq \frac29$. By proposition $2$, there are $p,q \in \mathbb Z^{\ast}$, $q>0$, $p$ and $q$ coprime, such that, if we set $A = A(p,q)$, $B = B(p,q)$ and $C = C(p,q)$, we have
$$X = \frac{A}{C}, \;\;\;\; Y = \frac{B}{C} \;\;\;\; \text{and} \;\;\;\; \frac{p}{q} = \frac{9Y - 2}{9X -1}.$$
$\bullet$ Now,
$$\frac{9Y - 2}{9X-1} = \frac{9\left(\frac{y}{z} \right) -2}{9 \left(\frac{x}{z} \right) -1} = \frac{9y - 2z}{9x - z} = I(S).$$
Hence, $I(S) = \frac{p}{q}$ (where $\frac{p}{q}$ is in lowest terms, with $q >0$). Next, let $\delta$ be the positive gcd of $A$ and $B$. By proposition $1$, $\delta$ divides $2 \cdot 3^5 \cdot 7 \cdot 59$, where we have $[ \ 7 \mid \delta \iff 7 \mid p \ ]$ and $[ \ 59 \mid \delta \iff 59 \mid q \ ]$.\\
Finally, we have $X = \frac{A}{C}$, hence $\frac{x}{z} = \frac{A/\delta}{C/\delta}$, where \textit{both} fractions are in lowest terms with $z>0$ and $\frac{C}{\delta}> 0$. We conclude that $x = \frac{A}{\delta}$ and $z = \frac{C}{\delta}$.\\
Similarly, we have $Y = \frac{B}{C}$, hence, $\frac{y}{z} = \frac{B/\delta}{C/\delta}$, where \textit{both} fractions are in lowest terms with $z > 0$ and $\frac{C}{\delta} > 0$. We conclude that $y = \frac{B}{\delta}$ (and $z = \frac{C}{\delta}$).
\end{proof}\hfill \\ \\
{\Large \bf{6. Proof of theorem 2:}} \\ \\
{\large \textbf{Lemma 6:}} \\
Let $(x,y,z)$ be a pp-solution of $(2)$, with $z \neq 9$ ($z > 0$). Note that the incidences of the $4$ primitive solutions of $(2)$, $(\pm x, \pm y, z)$, are all defined. Then, among these $4$ solutions, there is (at least) one, say $S = (x',y',z)$ (where $x' \in \{\pm x\}$ and $y' \in \{\pm y\}$), such that, if $I(S)$ is written in lowest terms as $I(S) = \frac{p}{q}$, with $q >0$, and, if $A = A(p,q)$, $B = B(p,q)$ and if $\delta$ is the positive gcd of $A$ and $B$, then, $\delta$ is neither divisible by $7$ nor by $59$.
\begin{proof} \hfill \\
We have $7x^2 + 59y^2 = 4z^2$, where $x$ and $y$ are positive and coprime (and $z >0$). \; Consider the four primitive solutions of $(2)$:
$$S_1 = (x,y,z), \;\;\; S_2 = (-x,y,z), \;\;\; S_3 = (x,-y,z) \;\;\; \text{and} \;\;\; S_4 = (-x,-y,z).$$
Since $z \neq 9$, then, each of these solutions is $\neq (\pm1,\pm2,9)$, so that its incidence is defined. Let us write $I(S_i)$, $i=1,2,3,4$, ${}$ in lowest terms as ${}$ $I(S_i) = \frac{p_i}{q_i}$, ${}$ with $q_i > 0$.\hspace{3.5mm}We have precisely \\ \\ \\

\[\setlength\arraycolsep{1pt}
\begin{array}[t]{rclcl@{\qquad}rclcl}
I(S_1) & = & \frac{9y-2z}{9x -z} & = & \frac{p_1}{q_1}, \,\,\,\,\,\,\,    &  I(S_2) & = & \frac{-9y + 2z}{9x + z} & = & \frac{p_2}{q_2} , \\ \\
I(S_3) & = & \frac{9y + 2z}{-9x + z}  & = & \frac{p_3}{q_3}, \,\,\,\,\,\,\,  & I(S_4) & = & \, \frac{9y + 2z}{9x + z} & = & \frac{p_4}{q_4}.
\end{array}
\]

Since the fractions $\frac{p_i}{q_i}$ are irreducible, there are $4$ integers $\lambda_i$ such that\\ \\
\null \hspace{2cm} $9y - 2z = \lambda_1p_1  \ \ \ \  (1,a)$ \hfill $9x - z = \lambda_1q_1 \ \ \ \ (1,b)$ \hspace{2cm} \\
\null \hspace{1.46cm} \; $-9y + 2z = \lambda_2p_2 \ \ \ \ (2,a)$ \hfill $9x + z = \lambda_2q_2 \ \ \ \ (2,b)$ \hspace{2cm} \\
\null \hspace{2cm} $9y + 2z = \lambda_3p_3 \ \ \ \ (3,a)$ \hfill $-9x + z = \lambda_3q_3 \ \ \ \ (3,b)$ \hspace{2cm} \\
\null \hspace{2cm} $9y + 2z = \lambda_4p_4 \ \ \ \ (4,a)$ \hfill $9x + z = \lambda_4q_4 \ \ \ \ (4,b)$ \hspace{2cm} \\ \\
For $i=1,2,3,4$, set $A_i = A(p_i,q_i)$, $B_i = B(p_i,q_i)$ and let $\delta_i$ be the positive gcd of $A_i$ and $B_i$. Recall (by proposition $1(ii)$) that, for $i=1,2,3,4$, we have:
$$[ \ 7 \mid \delta_i \iff 7\mid p_i \ ] \ \ \ \ \text{and} \ \ \ \ [ \ 59\mid \delta_i \iff 59 \mid q_i \ ]$$
Now, let us focus on the following arguments $(A1)$ and $(A2)$: \\ \\
\textbf{(A1)} \textit{If} $7$ were to divide $p_1$ and $p_3$, then, by $(1,a)$ and $(3,a)$, $7$ would divide $(9y - 2z) \pm (9y + 2z)$, that is, $7$ would divide $18y$ and $-4z$, hence, $7$ would divide $y$ and $z$, a contradiction since $\gcd(y,z) = 1$. Similarly, if $7$ were to divide $p_1$ and $p_4$, then using $(1,a)$ and $(4,a)$, we would obtain the same contradiction.\\
$\bullet$ Argument $(A1)$ shows that:
$$\textit{If \, \, $7 \mid p_1$, \, \, then, \, \, $7 \nmid p_3p_4$}$$ \hfill \\
\textbf{(A2)} \textit{If} $7$ were to divide $p_2$ and $p_3$, then, by $(2,a)$ and $(3,a)$, $7$ would divide $(-9y+2z) \pm (9y + 2z)$, that is, $7$ would divide $4z$ and $-18y$, hence, $7$ would divide $z$ and $y$, ${}$ a contradiction since $\gcd(y,z) = 1$. Similarly, \textit{if} $7$ were to divide $p_2$ and $p_4$, then using $(2,a)$ and $(4,a)$, we would obtain the same contradiction. \\
$\bullet$ Argument $(A2)$ shows that:
$$\textit{If \, \, $7 \mid p_2$, \, \, then, \, \, $7 \nmid p_3p_4$}$$ \hfill \\
From arguments $(A1)$ and $(A2)$, we deduce that:
$$\textit{Either \, \, $7 \nmid p_1 p_2$ \, \, or \, \, $7 \nmid p_3 p_4$}$$
Since [$ \, 7\mid p_i \iff 7\mid \delta_i \, $], we conclude the following:
\begin{equation}
\textbf{Either \, \, $(7 \nmid \delta_1 \ \ \textbf{and} \ \ 7 \nmid \delta_2)$ \ \ \ \textbf{or} \ \ \ $(7 \nmid \delta_3 \ \ \textbf{and} \ \ 7 \nmid \delta_4)$} \tag{\textbf{R1}}
\end{equation}
Next, we provide two similar arguments $(A3)$ and $(A4)$, relative to the prime $59$: \\ \\
\textbf{(A3)} \textit{If} $59$ were to divide $q_1$ and $q_2$, then, by $(1,b)$ and $(2,b)$, $59$ would divide $(9x - z) \pm (9x + z)$, that is, $59$ would divide $18x$ and $-2z$, hence, $59$ would divide $x$ and $z$, a contradiction since $\gcd(x,z) = 1$.\\
Similarly, \textit{if} $59$ were to divide $q_1$ and $q_4$, then, using $(1,b)$ and $(4,b)$, we would obtain the same contradiction.\\
$\bullet$ Argument $(A3)$ shows that:
$$\textit{If \, \, $59 \mid q_1$, \, \, then, \, \, $59 \nmid q_2q_4$}$$ \hfill \\
\textbf{(A4)} \textit{If} $59$ were to divide $q_3$ and $q_2$, then, by $(3,b)$ and $(2,b)$, $59$ would divide $(-9x + z) \pm (9x + z)$, that is, $59$ would divide $2z$ and $-18x$, hence, $59$ would divide $z$ and $x$, a contradiction, since $\gcd(x,z) = 1$. Similarly, \textit{if} $59$ were to divide $q_3$ and $q_4$, then using $(3,b)$ and $(4,b)$, we would obtain the same contradiction.\\
$\bullet$ Argument $(A4)$ shows that:
$$\textit{If \, \, $59 \mid q_3$, \, \, then, \, \, $59 \nmid q_2q_4$}$$
From arguments $(A3)$ and $(A4)$, we deduce that:
$$\textit{Either \, \, $59 \nmid q_1 q_3$ \, \, or \, \, $59 \nmid q_2q_4$}$$
Since [$\, 59 \mid q_i \iff 59 \mid \delta_i \, $], we conclude the following:
\begin{equation}
\textbf{Either \, \, $(59 \nmid \delta_1 \ \ \textbf{and} \ \ 59 \nmid \delta_3)$ \ \ \ \textbf{or} \ \ \ $(59 \nmid \delta_2 \ \ \textbf{and} \ \ 59 \nmid \delta_4)$} \tag{\textbf{R2}}
\end{equation}
Finally, by combining $(R1)$ and $(R_2)$ (that leads to $4$ cases), we see that for some $i \in \{1,2,3,4\}$, \\ $\delta_i$ is neither divisible by $7$ nor by $59$.
\end{proof} \hfill \\
Now we are in a position to prove the crucial result from which we will deduce theorem 2. \\ \\
{\large \textbf{Lemma 7:}}\\
Let $S_0 = (x,y,m)$ be a pp-solution of $(1)$, such that $S_0 \neq (1,2,5)$. Then, $S_0$ is one of the (two) successors of some $\textit{smaller}$ pp-solution of $(1)$.
\begin{proof}\hfill \\
Since $S_0 \neq (1,2,5)$, then, $m \ge 15$. By corollary $3$, $m$ is odd. Set $m = 2n + 1$ (so $n\ge 7$), and set $z = 3^n$ $(\ge 3^7)$. ${}$ From $7x^2 + 59y^2 = 3^{2n+1}$, we get
$$7x^2 + 59y^2 = 3z^2$$
As further $\gcd(x,y) = 1$, we see that $(x,y,z)$ is a pp-solution of $(2)$. Since $z \ge 3^7 > 9$, then, by lemma $6$, there are $x' \in \{\pm x\}$ and $y' \in \{\pm y\}$ such that, with $S = (x',y',z)$, if $I(S)$ is written in lowest terms as $I(S) = \frac{p}{q}$, with $q > 0$, if $A = A(p,q)$, $B = B(p,q)$, $C = C(p,q)$, and if $\delta$ is the positive gcd of $A$ and $B$, then $\delta$ is neither divisible by $7$ nor by $59$. But, by proposition $1$, $\delta$ divides $2 \cdot 3^5 \cdot 7 \cdot 59$. We conclude here that $\delta$ must divide $2 \cdot 3^5$. In such situation, we claim that $\delta$ cannot have the factor $2$: \textit{Indeed}, for the purpose of contradiction, suppose that $\delta = 2 \cdot 3^s$, $0 \le s \le 5$. By proposition $3$, we would get $z = \frac{C}{\delta}$. Hence, $C = \delta z = 2 \cdot 3^s \cdot 3^n = 2 \cdot 3^{n+s}$. That is, $9(7q^2 + 59p^2) = 2 \cdot 3^{n+s}$. Therefore, we would get
$$7q^2 + 59 p^2 = 2 \cdot 3^{n+s -2}$$
(where $n + s - 2 \ge 7 + 0 -2 = 5$).\\
However, this is a \textit{contradiction}, according to corollary $4$. \; Finally, $\delta$ must have the form\\ $\delta = 3^s$, ${}$ $0 \le s \le 5$.\; Recall by proposition $3$, that
$$x ' = \frac{A}{\delta}, \;\;\;\; y' = \frac{B}{\delta}, \;\;\;\; \text{and} \;\;\;\; z = \frac{C}{\delta}.$$
We consider two cases: \\ \\
\textit{Case 1}: $\delta = 1$\\
We then have ${}$ $x' = A(p,q)$,\; $y'=B(p,q)$ ${}$ and ${}$ $z = C(p,q)$.\\
From $z = 3^n = C(p,q) = 9(7q^2 + 59p^2)$, we obtain
$$7q^2 + 59p^2 = 3^{n-2}$$
Therefore, $(q,|p|, n-2)$ is a pp-solution of $(1)$. ${}$ Now, we have
$$x = |x'| = |A(p,q)| = |A(\pm|p|, q)|$$
and
$$y = |y'| = |B(p,q)| = |B(\pm|p|, q)|$$
Hence, in virtue of remark $3$ $(i)$, ${}$ we conclude that
$$(x,\ y,\ 2n+1) \text{ is the second successor of } (q, \ |p|, \ n-2).$$
(Note that $2n+1 >n-2$). \\ \\
\textit{Case 2}: $3 \mid \delta$\\
In particular, $3 \mid A$. \; If we had $pq \equiv 1 \pmod 3$, we would get
$$A = A(p,q) = 59p^2 - 236 pq - 7q^2 \equiv 59 - 236 - 7 = -184 \equiv -1 \pmod 3$$
which is in contradiction with $3 \mid A$. ${}$ Hence, $pq \equiv  -1 \pmod 3$. ${}$ By corollary $2$, we get $\delta = 3^5$. \\
Hence, ${}$ $x' = 3^{-5}A$, \; $y' = 3^{-5} B$ ${}$ and ${}$ $z = 3^{-5}C$. We may write
$$z = 3^n = 3^{-5}C = 3^{-5} C(p,q) = 3^{-5} \cdot 9(7q^2 + 59p^2).$$
Hence,
$$7q^2 + 59p^2 = 3^{n+3}$$
Therefore, $(q, \ |p|, \ n+3)$ is a pp-solution of $(1)$.\;\: Now, we have
$$x = |x'| = 3^{-5}|A(p,q)| = 3^{-5}|A(\pm|p|, q)|$$
and
$$y = |y'| = 3^{-5}|B(p,q)| = 3^{-5} |B(\pm|p|, q)|.$$
Hence, in virtue of remark $3$ $(ii)$, ${}$ we conclude that
$$(x, \ y, \ 2n+1) \text{ is the first successor of } (q, \ |p|, \ n+3).$$
(As $n \ge 7$, then, $2n + 1 > n+3$).
\end{proof} \hfill \\ \\
\textit{Proof of theorem 2}: \\
Recall that $S_1 = (1,2,5)$ is the \textit{smallest} pp-solution of $(1)$. \\
Let $T_1 = (x_1,y_1,m_1)$ be a pp-solution of $(1)$. If $T_1 = S_1$, there is nothing to prove. \textit{From now on}, we assume that $T_1 \neq (1,2,5)$. \\
By lemma $7$, ${}$ $T_1$ is a successor of a \textit{smaller} pp-solution of $(1)$, say
$$T_2 = (x_2,y_2, m_2), \, \, \, \, \text{where $m_1 > m_2$}.$$
If $T_2 = (1,2,5)$, we stop the process. Otherwise, by lemma $7$, ${}$ $T_2$ is a successor of a \textit{smaller} pp-solution of $(1)$, say
$$T_3 = (x_3,y_3,m_3), \,\,\,\, \text{where $m_2 > m_3$}.$$
This process must stop, otherwise we would obtain an infinite strictly decreasing sequence of positive integers $m_1 > m_2 > m_3 > \cdots $, which is impossible ${}$! \\
Hence, for some $r\ge 2$, $T_r = (x_r,y_r,m_r) = (1,2,5)$. In other words, $T_1$ is a successor of a successor $\ldots$ of a successor of $S_1$, \; so that ${}$ $T_1$ belongs to the binary tree described in corollary $1$. \hfill $\square$ \\ \\
{\large \textbf{Remark 5:}}\\
A direct consequence of theorem $2$ is the uniqueness of a pp-solution $(x,y,m_0)$ of $(1)$, for a \textit{given} suitable $m_0 = 10k_0 + 5$. \\ \\
$\bullet$ Finally, noting that in any solution $(x,y,m)$ of $(1)$, $x$ and $y$ have opposite parities we show that the case $[ \, x$ even and  $y$ odd$ \, ]$ never arises. \\ \\
{\large \textbf{Proposition 4:}}\\
Let $(x,y,m)$ be any solution (not necessarily primitive) of $(1)$. ${}$ Then, ${}$ $x$ is \textit{odd} and $y$ is \textit{even}.
\begin{proof} \hfill \\
Obviously, we only need to consider positive solutions of $(1)$. Further, since every positive solution of $(1)$ is obtained from a pp-solution of $(1)$ by multiplication by a power of $3$ (that is odd), it suffices then to prove the property for the pp-solutions of $(1)$:\\
By theorem $2$, the pp-solutions of $(1)$ lie all on the binary tree of the iterated successors of $S_1$.\\
We proceed by induction on the binary tree:\\
The property holds for $S_1$.\\
Next, we assume that the property holds for a pp-solution $S = (q,p,\omega)$ of $(1)$, and we prove that the property holds for the two successors of $S$. Let $S' = (x',y',m')$ and $S'' = (x'',y'',m'')$ denote respectively the first and second successor of $S$.\\
Since $B(p,q) = -118 p^2 - 14pq + 14q^2$ is even, and $y' = 3^{-5} |B(p,q)|$, then $y'$ is even (so $x'$ is odd). (Alternatively, from $q$ odd and $p$ even, we see that $A(p,q)$ is odd, so $x' = 3^{-5} |A(p,q)|$ is odd).\\
Similarly,\hspace{2mm}since\hspace{2mm}$B(-p,q) = -118p^2 + 14pq + 14q^2$ ${}$ is even, and $y'' = |B(-p,q)|$, then $y''$ is even \\ (so $x''$ is odd).
\end{proof} \hfill \\
{\large \textbf{Exercises:}}
\begin{enumerate}
\item[1.] The equation $3x^2 + 5y^2 = 7^m$ has no solution in integers $x,y,m$.
\item[2.] In any integer solution of $2x^2 + 3y^2 = 5^m$, ${}$ $m$ must be odd.
\item[3.] $(\star)$ Can one find distinct positive primes $a,b,c$ such that the equation $ax^2 + by^2 = c^z$ has exactly one primitive and positive solution?
\item[4.] $(\star)$ Let $(x,y,m)$ be a positive and primitive solution of $7x^2 + 59y^2 = 3^m$. Prove or disprove that $xy \equiv -1 \pmod3$.
\end{enumerate}
\hfill \\ \\ \\ \\ \\ \\ \\
{\Large \bf{References}} \\
\begin{enumerate}[label={[\arabic*]}]
\item W. Sierpinsky, \textit{On the equation} $3^x + 4^y = 5^z$, Wiadom. Mat. (2), 1955/1956, \ \textbf{1} \ 194-195 \\ (in Polish).
\item N. Terai, \textit{The Diophantine Equation} $a^x + b^y = c^z$, \ Proc. Japan Acad., 1994, \textbf{70A(2)} \ 22-26.
\item Bonyat Sroysang, \textit{More on the Diophantine Equation} $8^x + 19^y = z^2$, International Journal of Pure and Applied Mathematics, Vol.\hspace{1mm}81, \ No.\hspace{1mm}4, \ 2012, \: 601-604.
\item Julius Fergy T. Rabago, \textit{On the Diophantine Equation} $2^x + 17^y = z^2$, J. Indones. Math Soc., Vol.\hspace{1mm}22, No.\hspace{1mm}2 (2016) \ pp. 85-88.
\item Zhenfu Cao, Chuan I Chu and Wai Chee Shiu, \textit{The Exponential Diophantine Equation \\ $AX^2 + BY^2 = \lambda k^z$ and its Applications}, Taiwanese Journal of Mathematics, Vol.\hspace{1mm}12, No.\hspace{1mm}5, 2008, \ 1015-1034.

\item Borevitch Z.I., Shafarevich I.R., \textit{Number theory}, Acad. Press, 1966, Chap. 1, Par. 7. \\
\item Niven I., Zuckerman H.S., ${}$ \textit{An Introduction to the theory of Numbers}, ${}$ 4th Ed., ${}$ John Wiley \\ \& Sons.\\
\item Hardy G.H., Wright E.M., \textit{An Introduction to the Theory of Numbers}, Oxford, 5th Ed., 1979.\\
\item Tijdeman R., \textit{Exponential Diophantine Equations}, Proceedings of the International Congress of Mathematicians, Halsinki (381-387), 1978.\\
\item  Andreescu T., Andrica D., \textit{Number Theory, Structures, Examples, and Problems}, Birkhäuser, Boston-Basel-Berlin, 2009.\\
\item MORDELL, L.J., \textit{Diophantine Equations}, Acad. Press, London, 1969.\\
\item Samuel Pierre, \textit{Théorie Algébrique des Nombres}, Hermann, 1967.
\end{enumerate}
\end{document}